\documentclass[centertags,12pt]{amsart}
\usepackage{latexsym}
\usepackage{amsthm}
\usepackage{amssymb}
\usepackage{color}
\usepackage[utf8]{inputenc}
\usepackage{mathrsfs}
\usepackage[english]{babel}
\usepackage{setspace}
\usepackage{textcomp}
\usepackage{graphicx}
\usepackage{cite}
\usepackage{cancel}
\usepackage{comment}
\usepackage{hyperref}
\graphicspath{ {images/} }

\numberwithin{equation}{section}
\makeatletter
\renewcommand{\subsection}{\@startsection
{subsection}{2}{0mm}{\baselineskip}{-0.25cm}
{\normalfont\normalsize\bf}}
\makeatother

\newtheorem{theorem}{Theorem}[section]
\newtheorem{proposition}[theorem]{Proposition}
\newtheorem{lemma}[theorem]{Lemma}

\newtheorem*{theorem*}{Theorem}
   
\newtheorem{remark}[theorem]{Remark}

\newtheorem*{ThMain}{Main Theorem}

\begin{document}

\author[D. Bartoli]{Daniele Bartoli}
\address{
Dipartimento di Matematica e Informatica, Universit\`a degli Studi di Perugia}
	\email{daniele.bartoli@unipg.it}
\author[M. Timpanella]{Marco Timpanella}
\address{Dipartimento di Matematica, Informatica ed Economia, Universit\`a degli Studi della Basilicata}
	\email{marco.timpanella@unibas.it}

\title{A family of permutation trinomials in $\mathbb{F}_{q^2}$}

\begin{abstract}
Let $p>3$ and consider a prime power $q=p^h$. We completely characterize  permutation polynomials of $\mathbb{F}_{q^2}$ of the type $f_{a,b}(X) = X(1 + aX^{q(q-1)} + bX^{2(q-1)}) \in \mathbb{F}_{q^2}[X]$. In particular, using connections with algebraic curves over finite fields, we show that the already known sufficient conditions are also necessary. 
\end{abstract}

\maketitle

\section{Introduction}

Let $q=p^h$ be a prime power and denote by $\mathbb{F}_q$ the finite field with $q$ elements. A polynomial $f(x) \in  \mathbb{F}_q[X]$ is called a \emph{permutation polynomial} (PP) of $\mathbb{F}_q$ if it induces a permutation of $\mathbb{F}_q$. 

In this paper we consider   polynomials of the form
\begin{equation}\label{POL:family}
f_{a,b}(X) = X(1 + aX^{q(q-1)} + bX^{2(q-1)}) \in \mathbb{F}_{q^2}[X],
\end{equation}
where $a,b \in \mathbb{F}_{q^2}^*$.
This family of polynomials has been investigated by different authors; see for instance \cite{TZ2018,TZLH2018}.
Polynomials of type  \eqref{POL:family} belong to the more general family of PPs of $\mathbb{F}_{q^2}$ of the form
\begin{equation*}
F_{a,b,s_1,s_2}(X) = X(1 + aX^{s_1(q-1)} + bX^{s_2(q-1)}) \in \mathbb{F}_{q^2}[X],
\end{equation*}
where $1 \leq s_1, s_2 \leq q$. Such a family has been investigated in a number of papers, see for instance \cite{Bartoli2018,BG2018,BQ2018,DQWYY2015,GS2016,Hou2018,LP1997,LH2017,TZ2018,TZLH2018}. More in particular, necessary and sufficient conditions on $a, b$ for $F_{a,b,s_1,s_2}(X)$ to be a PP of $\mathbb{F}_{q^2}$ have been
determined only in the following cases:
\begin{enumerate}
    \item $(s_1, s_2) = (1, 2)$ and $p$ is arbitrary \cite{Hou2015};
    \item  $(s_1, s_2) = (-1/2, 1/2)$ and $p = 2$ \cite{TZ2018b};
    \item $(s_1, s_2) = (q, 2)$ and $p = 2$ \cite{Bartoli2018,HTZ2019,TZLH2018}.
\end{enumerate}

As mentioned above, in this paper we study polynomials $f_{a,b}(X)$ of the type \eqref{POL:family}. Know results are summarized in the following theorem.
\begin{theorem}
\begin{enumerate}
\item[(i)] Let $p=2$. Then $f_{a,b}(X)$ is PP of $\mathbb{F}_{q^2}$ if and only if 
$$b(1+a^{q+1}+b^{q+1})+a^{2q}=0$$
and 
$$\left\{
\begin{array}{ll}
Tr_{\mathbb{F}_q/\mathbb{F}_2}\left(1+1/a^{q+1}\right)=0& \textrm{ if } b^{q+1}=1\\
Tr_{\mathbb{F}_q/\mathbb{F}_2}\left(b^{q+1}/a^{q+1}\right)=0,& \textrm{ if } b^{q+1}\neq1,\\
\end{array}
\right.$$
see \cite{Bartoli2018,TZLH2018,Hou2018}.
\item [(ii)] Let $p=3$. Then $f_{a,b}(X)$ is PP of $\mathbb{F}_{q^2}$ if and only if 
$$\left\{
\begin{array}{ll}
a^qb^q=a(b^{q+1}-a^{q+1})\\
1-(b/a)^{q+1} \textrm{ is a square in } \mathbb{F}_q^*,\\
\end{array}
\right.$$
see \cite{TZ2018,HTZ2019}.
\item [(iii)] Let $p>3$. Then $f_{a,b}(X)$ is PP of $\mathbb{F}_{q^2}$ if either
\begin{equation}\label{PRIMA}
\left\{
\begin{array}{ll}
a^qb^q=a(b^{q+1}-a^{q+1})\\
1-4(b/a)^{q+1} \textrm{ is a square in } \mathbb{F}_q^*,\\
\end{array}
\right.
\end{equation}
or
\begin{equation}\label{SECONDA}
\left\{
\begin{array}{ll}
a^{q-1}+3b=0\\
-3(1-4(b/a)^{q+1}) \textrm{ is a square in } \mathbb{F}_q^*,\\
\end{array}
\right.
\end{equation}
see \cite{TZ2018}.
\end{enumerate}
\end{theorem}

The main aim of our paper is to show that for $p>3$ Conditions \eqref{PRIMA} and \eqref{SECONDA} are also necessary.

\begin{ThMain}

Let $p>3$ be a prime and $q=p^h$, with $h\geq 1$. The polynomial 
$$f_{a,b}(X) = X(1 + aX^{q(q-1)} + bX^{2(q-1)}) \in \mathbb{F}_{q^2}[X]$$
is a PP of $\mathbb{F}_{q^2}$ if and only if either Condition \eqref{PRIMA} or Condition \eqref{SECONDA} holds.
\end{ThMain}

To prove the main theorem we use a combination of a result contained in  \cite{PL2001,Zieve2009,Wang2007} (often called  AGW criterion) with (standard) techniques connected with algebraic curves \cite{Bartoli2018,BG2018}; see in particular Section \ref{Sec:GeneralIdea}.  

We would like to point out that another possible approach to establish whether a polynomial of type $f_{a,b}(X)$ permutes $\mathbb{F}_{q^2}$ is based on the classification of permutation rational functions of degree three on $\mathbb{P}_q^1=\mathbb{F}_q \cup \{\infty\}$; see \cite{FM2018}. Such a  technique, in principle, can be applied to study permutation properties of polynomials $x^rh(X^{q-1})\in \mathbb{F}_{q^2}[X]$,  using (if any) classifications of rational functions of a fixed degree. Although this method seems promising, it is not straightforward to obtain in this way conditions on $a,b\in \mathbb{F}_{q^2}$ for which the polynomial $f_{a,b}(X)$ permutes $\mathbb{F}_{q^2}$; see Section \ref{Sec:Appendix}.

\section{Necessary conditions connected with algebraic curves}\label{Sec:GeneralIdea}
It is well known that $f_{a,b}(X)$ is a PP of $\mathbb{F}_{q^2}$ if and only if $h_{a,b}(X) := X(1 + aX^q + bX^2)^{q-1}$ permutes the set $\mu_{q+1}$ of the $(q+1)$-roots of unity in $\mathbb{F}_{q^2}$; see for instance \cite{PL2001,Zieve2009,Wang2007}. Let $x\in \mu_{q+1}$ be such that $bx^3+x+a\neq 0$. Then $x^{q}=1/x$ and so 

$$h_{a,b}(x) := \frac{x(1 + ax^q + bx^2)^{q}}{1+ax^q+bx^2}=\frac{a^qx^3+x^2+b^q}{bx^3+x+a}.$$

Consider now the function 
\begin{equation}\label{Eq:g}
g_{a,b}(X)=\frac{N_{a,b}(x)}{D_{a,b}(x)}=\frac{a^qX^3+X^2+b^q}{bX^3+X+a}\in \mathbb{F}_{q^2}(X).
\end{equation} Recall that the main aim of the present paper is  to obtain necessary conditions for $f_{a,b}(X)$ to be PP in $\mathbb{F}_{q^2}$. If $f_{a,b}(X)$ permutes $\mathbb{F}_{q^2}$ then $h_{a,b}(X)$ permutes $\mu_{q+1}$ and so does $g_{a,b}(X)$ (also, $bX^3+X+a$ has no roots in $\mu_{q+1}$).

In this paper we follow an approach based on the study of algebraic curves associated with $f_{a,b}(X)$. Consider the quartic curve $\mathcal{C}_{a,b}$ with the affine equation 
\begin{equation}\label{Eq:F}
F_{a,b}(X,Y)=\frac{(a^qX^3+X^2+b^q)(bY^3+Y+a)-(a^qY^3+Y^2+b^q)(bX^3+X+a)}{X-Y}=0.
\end{equation}
This curve is defined over $\mathbb{F}_{q^2}$. 
Note that $g_{a,b}(X)$ permutes $\mu_{q+1}$ if and only if there is no $(\overline{u},\overline{v})\in \mu_{q+1}^2$, $\overline{u}\neq \overline{v}$ be such that $g_{a,b}(\overline{u})=g_{a,b}(\overline{v})$.

Let $e\in \mathbb{F}_{q^2}\setminus \mathbb{F}_q$ such that $e^q=-e$ and  denote by $\psi$ and $\phi$ the transformations 
$$\psi(X,Y)=\left(\frac{X+e}{X-e},\frac{Y+e}{Y-e}\right), \qquad 
\phi(X,Y)=\left(e\frac{X-1}{X+1},e\frac{Y-1}{Y+1}\right).$$

Now, $G_{a,b}(X,Y)=(X-e)^2(Y-e)^2F_{a,b}(\psi(X,Y))\in \mathbb{F}_{q}[X,Y]$ and it defines a  curve $\mathcal{D}_{a,b}$. Since $(X-1)^2(Y-1)^2G_{a,b}(\phi(X,Y)) = 16e^4F_{a,b}(X,Y)$, the curves $\mathcal{C}_{a,b}$ and $\mathcal{D}_{a,b}$ are $\mathbb{F}_{q^2}$-isomorphic.

\begin{lemma}\label{Lemma}
Let $q\geq 47$. If $f_{a,b}(x)\in \mathbb{F}_{q^2}[x]$ is a PP then $\mathcal{C}_{a,b}$ is not absolutely irreducible. 
\end{lemma}
\proof
If $\mathcal{C}_{a,b}$ is  absolutely irreducible then  $\mathcal{D}_{a,b}$  is absolutely irreducible too. Since $\mathcal{D}_{a,b}$ is defined over $\mathbb{F}_q$ and it has degree at most $4$, by  Hasse-Weil Theorem (see \cite[Theorem 5.2.3]{Sti})), $\mathcal{D}_{a,b}$ possesses affine $\mathbb{F}_q$-rational points $(x_0,y_0)$ with $x_0\neq y_0$ whenever 
\begin{equation}\label{Eq:HW}
q+1-6\sqrt{q}-6=q-6\sqrt{q}-5>0,
\end{equation}
where $6$ is the maximum number of points of $\mathcal{D}_{a,b}$ lying on $X=Y$ or on the line at infinity (in fact, $\mathcal{D}_{a,b}$  has only two ideal points). Equation \eqref{Eq:HW} is satisfied if  $q\geq 47$. So there exists an $\mathbb{F}_{q}$-rational point $(x_0,y_0)$, with $x_0\neq y_0$, of $\mathcal{D}_{a,b}$. Such a point  gives rise to a pair $\left(\frac{x_0+e}{x_0-e},\frac{y_0+e}{y_0-e}\right)\in \mu_{q+1}^2$ with $ \frac{x_0+e}{x_0-e}\neq \frac{y_0+e}{y_0-e}$ satisfying 
$$g_{a,b}\left(\frac{x_0+e}{x_0-e}\right)=g_{a,b}\left(\frac{y_0+e}{y_0-e}\right).$$
This is a contradiction to $f_{a,b}(x)$ being a PP in $\mathbb{F}_{q^2}$.
\endproof

Due to Lemma \ref{Lemma}, a necessary condition for $f_{a,b}(X)$ to be PP is that, for $q\geq 47$, $\mathcal{C}_{a,b}$ is reducible. 

\begin{remark}
An exhaustive computer search (using the software MAGMA \cite{MAGMA}) shows that for $q\in \{5,7,11,13,17,19,23,25,29,31,37,41,43\}$ Conditions \eqref{PRIMA} and \eqref{SECONDA} are necessary for a polynomial $f_{a,b}(X)$ to be a PP of $\mathbb{F}_{q^2}$.
\end{remark}

In the next Sections, we distinguish the cases depending on the degree of $\mathcal{C}_{a,b}$, which is connected with the degree of $GCD(a^qX^3+X^2+b^q,bX^3+X+a)$; see in particular Propositions \ref{Prop:GCD2} and \ref{Prop:GCD1} for the cases $\deg(GCD(a^qX^3+X^2+b^q,bX^3+X+a))\in \{1,2\}$. In the case of $\deg(GCD(a^qX^3+X^2+b^q,bX^3+X+a))=0$ we consider all the possibilities for the factors of $\mathcal{C}_{a,b}$; see Theorems \ref{Th:rette}, \ref{Th:coniche1}, and  \ref{Th:coniche2}.

To conclude this section, we include here an equivalent description of Conditions \eqref{PRIMA} and \eqref{SECONDA} which will be helpful in the rest of the paper. As a notation, the symbol $\square_q^*$ denotes the set of all nonzero squares in $\mathbb{F}_q$.

\begin{proposition}
Let $p>3$. If
\begin{equation}\label{PRIMABis}
a \in \mathbb{F}_{q^2}^*, \qquad b=v/a^2,\qquad v\in \mathbb{F}_q^* \qquad v^2-a^{q+1}v-a^{3q+3}=0, \qquad -3a^{2q+2}-4v\in \square_q^*
\end{equation}
then Condition \eqref{PRIMA} holds.
If 
\begin{equation}\label{SECONDABis}
b=-a^{q-1}/3, \qquad 3a^{q+1}(4-9a^{q+1})\in \square_q^*
\end{equation}
then Condition \eqref{SECONDA} holds.
\end{proposition}
\proof
First, we prove that Condition \eqref{PRIMABis} yields Condition \eqref{PRIMA}. Indeed, it is readily seen that if $b=v/a^2$ and $v^2-a^{q+1}v-a^{3q+3}=0$, then $a^{q-1}b^q=b^{q+1}-a^{q+1}$ holds. Also, $b=v/a^2$ and  $v^2-a^{q+1}v-a^{3q+3}=0$ yield  $$1-4(b/a)^{q+1}=\frac{a^{3q+3}-4v^2}{a^{3q+3}}=\frac{-3a^{2q+2}-4v}{a^{2q+2}}.$$ Therefore, since $-3a^{2q+2}-4v\in \square_q^*$, also $1-4(b/a)^{q+1}$ is a square in $\mathbb{F}_q^*$.

Finally, we show that Condition \eqref{SECONDABis} yields Condition \eqref{SECONDA}. Indeed, $b=-a^{q-1}/3$ yields $b^{q+1}=1/9$ and we have 
$$-3(1-4(b/a)^{q+1})=-3(1-\frac{4}{9a^{q+1}})=\frac{3a^{q+1}(4-9a^{q+1})}{9a^{2q+2}}.$$  Therefore, since $3a^{q+1}(4-9a^{q+1})\in \square_q^*$, also $-3(1-4(b/a)^{q+1})$ is a square in $\mathbb{F}_q^*$.

\endproof

Actually, it is possible to prove that Conditions \eqref{PRIMA} and \eqref{SECONDA} are equivalent to Conditions \eqref{PRIMABis} and \eqref{SECONDABis}: we do not need this implication in our argument.

\section{CASE \texorpdfstring{ $\deg(GCD(a^qX^3+X^2+b^q,bX^3+X+a))=2$}{Lg}}

Assume that $\deg(GCD(a^qX^3+X^2+b^q,bX^3+X+a))=2$. This means that 
\begin{eqnarray*}
GCD(a^qX^3+X^2+b^q,bX^3+X+a)=GCD(a^qX^3+X^2+b^q,bX^2 - a^qX - a^{q+1} + b^{q+1})\\
=GCD((a^{2q}+b)X^2 +( a^{2q+1} - a^qb^{q+1})X + b^{q+1},bX^2 - a^qX - a^{q+1} + b^{q+1}).
\end{eqnarray*}
Since the GCD must have degree two, these two last polynomial must be proportional, that is 
\begin{equation}\label{Condizione_GCD2}
a^q(-a^{q+1}b - a^{2q} +  b^{q+2} - b)=0,  \quad 
a^{2q}(-a^{q+1} - b/a^{q-1} + b^{q+1})=0.
\end{equation}

Since $ab\neq 0$, from the last condition we get $b^{q+1}-a^{q+1}=b/a^{q-1} \in \mathbb{F}_q$ and then $b/a^{q-1}=b^qa^{q-1}$, that is $b^{q-1}=(1/a^{2})^{q-1}$. This is equivalent to $b=v/a^2$ with $v\in \mathbb{F}_{q}^*$. Also, Conditions \eqref{Condizione_GCD2} yields $v^2-a^{q+1}v-a^{3q+3}=0$. A direct check shows that $a \in \mathbb{F}_{q^2}^*$, $b=v/a^2$, $v^2-a^{q+1}v-a^{3q+3}=0$ satisfy \eqref{Condizione_GCD2}.

\begin{proposition}\label{Prop:GCD2}
If $f_{a,b}$ is a PP and  $\deg(GCD(a^qX^3+X^2+b^q,bX^3+X+a))=2$ then \eqref{PRIMABis} holds.
\end{proposition}
\proof
We already saw that $\deg(GCD(a^qX^3+X^2+b^q,bX^3+X+a))=2$ is equivalent to 
$a \in \mathbb{F}_{q^2}^*$, $b=v/a^2$,  $v^2-a^{q+1}v-a^{3q+3}=0$. The common factor between $a^qX^3+X^2+b^q$ and $bX^3+X+a$ is $r_1=bX^2 - a^qX - a^{q+1} + b^{q+1}$. The roots of $r_1$ are 
$$Z=\frac{a^q+\delta}{2b}, \quad \delta^2=4a^{q+1}b + a^{2q} - 4b^{q+2}=\frac{-3a^{4q+4} - 4a^{2q+2}v}{a^{2q+4}}=\frac{-3a^{2q+2}-4v}{a^2}.$$
This means that $\delta=\sqrt{-3a^{2q+2}-4v} /a$. Let $\eta=\sqrt{-3a^{2q+2}-4v}$. Suppose that  $-3a^{2q+2}-4v \notin \square_q^*$. Then $\delta^q=-\eta /a^q$ and $Z^{q}= \frac{a-\eta/a^q}{2b^q} $. So 
\begin{eqnarray*}
Z^q=1/Z &\iff& \frac{a-\eta/a^q}{2b^q}=\frac{2b}{a^q+\eta/a} \iff a^{q+1}-\eta^2 /a^{q+1}=4b^{q+1}\\
&\iff&  a^{q+1}-\eta^2 /a^{q+1}-4v^2/a^{2q+2}=0\\
&\iff& a^{3q+3}-(-3a^{2q+2}-4v) a^{q+1}-4v^2=0\\
&\iff& 4a^{3q+3}+4va^{q+1}-4v^2=0,\\
\end{eqnarray*}
and therefore the roots of $r_1$ are in $\mu_{q+1}$. This is a contradiction to $f_{a,b}$ permuting $\mathbb{F}_{q^2}$. 

\endproof

\section{CASE \texorpdfstring{ $\deg(GCD(a^qX^3+X^2+b^q,bX^3+X+a))=1$}{Lg}}

\begin{proposition}\label{Prop:GCD1}
If $f_{a,b}$ is a PP and  $\deg(GCD(a^qX^3+X^2+b^q,bX^3+X+a))=1$ then $f_{a,b}(X)$ is not a PP. 
\end{proposition}
\proof
 We already saw in the previous subsection that 
$(GCD(a^qX^3+X^2+b^q,bX^3+X+a))=GCD((a^{2q}+b)X^2 +( a^{2q+1} - a^qb^{q+1})X + b^{q+1},bX^2 - a^qX - a^{q+1} + b^{q+1})$. 

We distinguish two cases.
\begin{itemize}
\item If $b=-a^{2q}$ the GCD is $(a^{2q}+b)X^2 +( a^{2q+1} - a^qb^{q+1})X + b^{q+1}=a^{2q+1}((a^{q+1}-1)X - a)$. In particular $a^{q+1}\neq 1$ and its root is $a/(a^{q+1}-1)$. Also, $bX^2 - a^qX - a^{q+1} + b^{q+1}=a^q(X + a)(a^q X - a^{q+1} + 1)$. So either $a/(a^{q+1}-1)=-a$ or $a/(a^{q+1}-1)=(a^{q+1}-1)/a^q$. In the former case $a=0$, a contradiction. In the latter case $a/(a^{q+1}-1)\in \mu_{q+1}$ a contradiction to $f_{a,b}(x)$ permuting $\mathbb{F}_{q^2}$.
\item If $b\neq-a^{2q}$ then the GCD must be $\ell(X)=(a^{2q+1}b + a^{3q} - a^qb^{q+2} + a^qb)X + a^{3q+1} + a^{q+1}b - a^{2q}b^{q+1}$. Since it must also be a factor of $bX^2 - a^qX - a^{q+1} + b^{q+1}$, 
$$ba^{2q}(a^{2q}+ b)(a^{3q+3} - 3a^{2q+2}b^{q+1} - a^{2q+2} - a^2 b + 3a^{q+1}b^{2q+2} - a^{q+1}b^{q+1} -   a^{2q}b^q - b^{3q+3} + 2b^{2q+2} - b^{q+1})=0.$$

So, the last factor must vanish. Consider the root $Z= -\frac{a^{2q+1}b + a^{3q} - a^qb^{q+2} + a^qb}{a^{3q+1} + a^{q+1}b - a^{2q}b^{q+1}}$ of $\ell(X)$. Since 
$$a^{3q+3} - 3a^{2q+2}b^{q+1} - a^{2q+2} - a^2 b + 3a^{q+1}b^{2q+2} - a^{q+1}b^{q+1} -   a^{2q}b^q - b^{3q+3} + 2b^{2q+2} - b^{q+1}=0,$$
$Z$ belongs to $\mu_{q+1}$. This is clearly not possible since no roots of $a^qX^3+X^2+b^q$ can belong to $\mu_{q+1}$, otherwise $f_{a,b}(x)$ does not permute $\mathbb{F}_{q^2}$.

\end{itemize}
\endproof

\section{CASE \texorpdfstring{ $\deg(GCD(a^qX^3+X^2+b^q,bX^3+X+a))=0$}{Lg}}

Recall that $\mathcal{C}_{a,b}$ and $\mathcal{D}_{a,b}$ are $\mathbb{F}_{q^2}$-isomorphic and therefore they have the same number of components. If $f_{a,b}(X)$ is a PP of $\mathbb{F}_{q^2}$ and $q\geq 47$ then $\mathcal{D}_{a,b}$ must have no affine $\mathbb{F}_q$-rational points off $X=Y$ (see Lemma \ref{Lemma}). In the following subsections we distinguish two cases depending on the degree of $\mathcal{D}_{a,b}$. 

\subsection{SUBCASE \texorpdfstring{$ \deg(\mathcal{D}_{a,b})=4$}{Lg}}

In this case the curve $\mathcal{D}_{a,b}$ must split either into four lines or two conics, otherwise at least one of its absolutely irreducible components is defined  over $\mathbb{F}_q$ and, by Hasse-Weil Theorem it contains affine $\mathbb{F}_q$-rational points off the line $X=Y$. 

Therefore $\mathcal{D}_{a,b}$ possesses either two or four absolutely irreducible components. So $\mathcal{C}_{a,b}$ splits into either two or four components, that is either four lines or two conics. All these cases are treated in Theorems \ref{Th:rette}, \ref{Th:coniche1}, and \ref{Th:coniche2}. Finally, since $(X,Y)\mapsto (Y,X)$ is an automorphism of $\mathcal{C}_{a,b}$, the unique cases to be considered are the following. 
\begin{itemize}
\item Four lines. In this case 
$$F_{a,b}(X,Y)=-b(X+A)(X+B)(Y+A)(Y+B),$$
for some $A,B\in \overline{\mathbb{F}_q}$; see Theorem \ref{Th:rette}.
\item Two conics switched by  $(X,Y)\mapsto (Y,X)$. So 
$$F_{a,b}(X,Y)=-b(X^2+AX+BY+C)(Y^2+AY+BX+C),$$ for some $A,B,C\in \overline{F}_q$; see Theorem \ref{Th:coniche0}.
 \item Two conics fixed by  $(X,Y)\mapsto (Y,X)$. So 
$$F_{a,b}(X,Y)=-b(XY+A(X+Y)+C)(XY+B(X+Y)+D),$$ for some $A,B,C,D\in \overline{F}_q$; see Theorem \ref{Th:coniche1}.
 \item Two conics switched by  $(X,Y)\mapsto (Y,X)$. So 
$$F_{a,b}(X,Y)=-b(XY+AX+BY+C)(XY+BX+AY+C),$$ for some $A,B,C\in \overline{F}_q$; see Theorem \ref{Th:coniche2}. 

\end{itemize}

\begin{theorem}\label{Th:rette}
Let 
$F_{a,b}(X,Y)$ as in \eqref{Eq:F}. If $F_{a,b}(X,Y)$ factorizes into four lines then $N_{a,b}(x)$ and $D_{a,b}(x)$ in  \eqref{Eq:g} share a factor.

\end{theorem}
\proof
Suppose that 
$$F_{a,b}(X,Y)=-b(X+A)(X+B)(Y+A)(Y+B),$$
for some $A,B\in \overline{F}_q$. Then 
\begin{eqnarray*}
a^{q+1} - b^{q+1} + bAB=0,\\
a + bA^2B + bAB^2=0,\\
a^q + bA + bB=0,\\
bA^2B^2 - b^q=0,\\
a^{q+1} - b^{q+1} + bA^2 + 2bAB + bB^2 + 1=0.
\end{eqnarray*}
Using $B=-(a^q + bA)/b$ one gets

\begin{eqnarray}
b(a^{q+1}b + a^{2q} - b^{q+2} + b)=0,\label{Prima}\\
b(a^{2q}A^2 + 2 a^q b A^3 + b^2A^4 - b^{q+1})=0,\nonumber\\
b(a^{q+1} - a^qA - b^{q+1} - bA^2)=0,\nonumber\\
b(ab + a^{2q}A + a^qbA^2)=0\nonumber,
\end{eqnarray}
so 
\begin{equation}\label{Seconda}
ab + a^{2q}A + a^qbA^2=0.
\end{equation}
 Also,  $A=-(a^q + bB)/b$  yields 
$ab + a^{2q}B + a^qbB^2=0$, so both $A$ and $B$ are roots of 
$$ a^qbT^2+a^{2q}T+ab.$$
Now, using System \eqref{Prima} and Equation \eqref{Seconda}, 
\begin{eqnarray*}
b^4(a^{2q+1} + ab - a^qb^{q+1})^2=0,\\
b^8(a^2b - a^{2q}b^q)^2=0,\\
b^2(a^{q+1}b + a^{2q} - b^{q+2} + b)^2=0.\\
\end{eqnarray*}
Condition $a^2b - a^{2q}b^q=0$ is equivalent to $b=v/a^2$, with $v\in \mathbb{F}_q^*$. Then, from $a^{2q+1} + ab - a^qb^{q+1}=0$,
$$a^{3q+3} + a^{q+1}v - v^2=0.$$
Now, the resultant between $N_{a,b}(x)$ and $D_{a,b}(x)$ in \eqref{Eq:g} with respect to $x$ reads
$$(a^{3q+3} + a^{q+1}v - v^2)^2(a^{3q+3} - a^{2q+2} - 2a^{q+1}v - v^2)\equiv 0,$$
and so $N_{a,b}(x)$ and $D_{a,b}(x)$ share a factor, a contradiction to $\deg(GCD(a^qX^3+X^2+b^q,bX^3+X+a))=0$. 
\endproof

\begin{theorem}\label{Th:coniche0}
Let 
$F_{a,b}(X,Y)$ as in \eqref{Eq:F}. 
If $F_{a,b}(X,Y)$ factorizes as
$-b(X^2+AX+BY+C)(Y^2+AY+BX+C),$ for some $A,B,C\in \overline{F}_q$, then $F_{a,b}(X,Y)$ factorizes into four linear factors. 
\end{theorem}
\proof From $F_{a,b}(X,Y)=-b(X^2+AX+BY+C)(Y^2+AY+BX+C)$ one gets
\begin{eqnarray*}
bC^2 - b^q=0\\
a + bAC + bBC=0\\
a^{q+1} - b^{q+1} + bAB + bC=0\\
a^{q+1} - b^{q+1} + bA^2 + bB^2 + 1=0\\
Bb=0\\
a^q + bA=0,
\end{eqnarray*}
and hence $B=0$. Since the two conics $X^2+AX+C$ and 
$Y^2+AY+C$ factorize in the product of two lines, the claim follows.
\endproof

\begin{theorem}\label{Th:coniche1}
Let 
$F_{a,b}(X,Y)$ as in \eqref{Eq:F}. 
If $F_{a,b}(X,Y)$ factorizes as
$-b(XY+A(X+Y)+C)(XY+B(X+Y)+D)$, for some $A,B,C,D\in \overline{F}_q$, then $N_{a,b}(x)$ and $D_{a,b}(x)$ share a factor.

\end{theorem}
\proof
From $F_{a,b}(X,Y)=-b(XY+A(X+Y)+C)(XY+B(X+Y)+D)$ one gets
\begin{eqnarray*}
a^{q+1} - b^{q+1} + bAB=0\\
C D - b^{q-1}=0\\
a^q + b A + bB=0\\
a + b A D + b B C=0\\
a^{q+1} - b^{q+1}+ 2bAB + bC + bD + 1=0.
\end{eqnarray*}
Since $B=-(a^q+bA)/b$,
\begin{eqnarray*}
CD - b^{q-1}=0\\
b(a - a^qC - bAC + bAD)=0\\
b(a^{q+1} - a^qA - b^{q+1} - bA^2)=0\\
b(a^{q+1} - 2a^qA - b^{q+1} - 2bA^2 + bC + bD + 1)=0.
\end{eqnarray*}
Both $C$ and $D$ are different from $0$ otherwise $b=0$. We use now $D=b^{q-1}/C$,
\begin{eqnarray*}
b^2(aC - a^qC^2 - bAC^2 + b^qA)=0\\
b^2(a^{q+1}C - 2a^qAC - b^{q+1}C - 2bA^2C + bC^2 + b^q + C)=0\\
b(a^{q+1}- a^qA - b^{q+1} - bA^2)=0.
\end{eqnarray*}
Now we consider $a^{q+1} - a^qA - b^{q+1} - bA^2=0$. Combining it with the other conditions in the above system (eliminating $A$) one gets
\begin{eqnarray*}
b^4( (- a^{q+1}b^2 + b^{q+3})C^4 
- a^{q+1}bC^3
+ (2a^{q+1}b^{q+1} +a^2b  + a^{2q}b^q  - 2b^{2q+2})C^2\\  
- a^{q+1}b^qC 
+ b^{3q+1}- a^{q+1}b^{2q}
)^2=0\\
b^6(a^{q+1}C - b^{q+1}C - bC^2 - b^q - C)^2=0.
\end{eqnarray*}
Finally, the resultant between these two last equations with respect to $C$ is $h(a,b)=$
$$b^{2q+10}(a^{3q+3} - 3a^{2q+2}b^{q+1} - a^{2q+2} - a^2b + 3a^{q+1}b^{2q+2} - 
            a^{q+1}b^{q+1} - a^{2q}b^q - b^{3q+3} + 2b^{2q+2} - b^{q+1})^2.$$
Note that the resultant between $N_{a,b}(x)$ and $D_{a,b}(x)$ with respect to $x$ is precisely $h(a,b)$. This means that if $F_{a,b}(X,Y)=-b(XY+A(X+Y)+C)(XY+B(X+Y)+D)$ then $N_{a,b}(x)$ and $D_{a,b}(x)$ share a factor.
\endproof

\begin{theorem}\label{Th:coniche2}
Let 
$F_{a,b}(X,Y)$ as in \eqref{Eq:F}.
If $F_{a,b}(X,Y)$ factorizes as
$-b(XY+AX+BY+C)(XY+BX+AY+C)$, for some $A,B,C\in \overline{F}_q$, then \eqref{SECONDABis} holds.
\end{theorem}
\proof
From $F_{a,b}(X,Y)=-b(XY+AX+BY+C)(XY+BX+AY+C)$ one gets
\begin{eqnarray*}
a^{q+1} - b^{q+1} + bAB=0\\
a + bAC + bBC=0\\
a^q + bA + bB=0\\
bC^2 - b^q=0\\
a^{q+1} - b^{q+1} + bA^2 + bB^2 + 2bC + 1=0.
\end{eqnarray*}
Since $B=-(a^q+bA)/b$,
\begin{eqnarray*}
b(a - a^qC)=0\\
b(a^{q+1}b + a^{2q} + 2a^qbA - b^{q+2} + 2b^2A^2 + 2b^2C + b)=0\\
b(a^{q+1} - a^qA - b^{q+1} - bA^2)=0\\
b(bC^2 - b^q)=0.
\end{eqnarray*}
From $C=1/a^{q-1}$,
\begin{eqnarray*}
a^2b - a^{2q}b^q=0\\
b(a^{q+1} - a^qA - b^{q+1} - bA^2)=0\\
b(a^{2q+1}b + 2ab^2 + a^{3q} + 2a^{2q}bA - 
            a^qb^{q+2}+ 2a^qb^2A^2 + a^qb)=0.
\end{eqnarray*}

From $a^{q+1} - a^qA - b^{q+1}- bA^2=0$,
\begin{eqnarray*}
(a^2b - a^{2q}b^q)^2=0\\  
b^4(3a^{2q+1}b + 2ab^2 + a^{3q} - 3a^qb^{q+2}
            + a^qb)^2=0.
\end{eqnarray*}
By $b^{q-1}=1/a^{2(q-1)}$, $b=v/a^2$ for some $v \in \mathbb{F}_q^*$, and then 
$$v^4(a^{q+1}+3v)^2(a^{3q+3}+a^{q+1}v-v^2)=0.$$
As we already saw, $a^{3q+3}+a^{q+1}v-v^2=0$ yields that $N_{a,b}(x)$ and $D_{a,b}(x)$ share a factor.
So, $a^{q+1}+3v=0$, that is $b=-a^{q-1}/3$,
\begin{equation}
3a^{q}A^2-9a^{q+1}A+9a^{q+2}-a=0, \qquad B=3a-A, \qquad C=1/a^{q-1}.
\end{equation}
In particular, since the discriminant of the polynomial $3a^{q}A^2-9a^{q+1}A+9a^{q+2}-a$ with respect to $A$ is  $\delta=3a^{q+1}(4-9a^{q+1})$ and it is a square in $\mathbb{F}_{q^2}$ (it belongs to $\mathbb{F}_q$),  $A,B,C \in \mathbb{F}_{q^2}$.

Since $f_{a,b}(x)$ is a $PP$ of $\mathbb{F}_{q^2}$, the  $\mathbb{F}_{q^2}$-rational conic $\mathcal{Q}$ of the type 
$$XY+AX+(3a-A)Y+1/a^{q-1}=0, \qquad \textrm{with } \quad 3a^{q}A^2-9a^{q+1}A+9a^{q+2}-a=0,$$
must not contain a point $P=(\alpha,\beta) \in \mu_{q+1}^2$, with $\alpha\neq \beta$. 

Suppose now that $\delta\in \mathbb{F}_q\setminus \square_q^*$, that is $\delta=\eta^2$, with $\eta\in \mathbb{F}_{q^2}\setminus\mathbb{F}_q^*$ and $\eta^q=-\eta$. Then 
$$A=\frac{9a^{q+1}+\eta}{6a^q}, \qquad A^q=\frac{9a^{q+1}-\eta}{6a}.$$
Let $\alpha\in \mu_{q+1}$. The corresponding $\beta\in \mu_{q+1}$ is 
$$\beta=\frac{-A\alpha-1/a^{q-1}}{\alpha+(3a-A)}=\frac{-6a^q\alpha - 9a^{q+1} + \eta}{9a^{q+1}\alpha + \eta \alpha + 6a}.$$
By direct checking, 
$$\beta^{q+1}=\frac{(-6a/\alpha - 9a^{q+1} - \eta)(-6a^q\alpha - 9a^{q+1} + \eta)}{(9a^{q+1}\alpha - \eta /\alpha + 6a^q)(9a^{q+1}\alpha + \eta \alpha + 6a)}=1,$$
and so $\beta\in \mu_{q+1}$. Since at most two points $(x,x)$, with $x\in \mu_{q+1}$, belong to $\mathcal{Q}$, whenever $q\geq3$ there exists at least one point $(\alpha,\beta)\in (\mu_{q+1})^2$, with $\alpha\neq \beta$.
So,  $3a^{q+1}(4-9a^{q+1})\in \square_q^*$.
\endproof

\subsection{SUBCASE \texorpdfstring{$ \deg(\mathcal{D}_{a,b})<4$}{Lg}}
By direct computation, 
\begin{eqnarray*}
G_{a,b}(X,Y)&=&(3a^{q+1} + 2a + 2a^q - 3b^{q+1} - b - b^q + 1)X^2Y^2\\
&&+(-2a e + 2a^q e - 2b e + 2 b^q e)(X^2 Y+XY^2)\\
&&+(a^{q+1} e^2 - b^{q+1} e^2 - b e^2 - b^q e^2 - e^2)(X^2+Y^2)\\
&&+(-8a^{q+1} e^2 + 8 b^{q+1} e^2 - 4b e^2 - 4b^q e^2)XY\\
&&(2ae^3 - 2a^qe^3 - 2be^3 + 2b^qe^3)(X+Y)\\
&&+(3a^{q+1}e^4 - 2ae^4 - 2a^qe^4 - 3b^{q+1}e^4 - be^4 - b^qe^4 + e^4).
\end{eqnarray*}
If $3a^{q+1} + 2a + 2a^q - 3b^{q+1} - b - b^q + 1=0$ and $-2a e + 2a^q e - 2b e + 2 b^q e\neq 0$  then $G_{a,b}(X,Y)$ has degree $3$ and it contains an absolutely irreducible components defined over $\mathbb{F}_q$. To see this, it is enough to observe that $(1,0,0)$ is  simple $\mathbb{F}_q$-rational point of $\mathcal{D}_{a,b}$. So, if $\mathcal{D}_{a,b}$ has degree less than $4$ and it has no absolutely irreducible $\mathbb{F}_q$-rational component, then   
the terms of $G_{a,b}(X,Y)$ of degree 3 and 4 must vanish, that is
\begin{eqnarray*}
3a^{q+1} + 2a + 2a^q - 3b^{q+1} - b - b^q + 1=0\\
-a  + a^q  - b  + b^q =0,
\end{eqnarray*}
whence 
\begin{equation*}\label{condizConica}
(3a^q - 3b + 1)(a + b + 1)=0.
\end{equation*}

Assume first $a=-(b+1)$. In this case $f_{a,b}(1)=0$ and $f_{a,b}$ is not a permutation, a contradiction.
Then $b=a^q+1/3$ and $\mathcal{D}_{a,b}$ is the curve with affine equation $G_{a,b}^{\prime}(X,Y)=0$, where
\begin{eqnarray*}
G_{a,b}^{\prime}(X,Y)&=&(3a + 3a^q + 4)(X^2+Y^2)\\
&&+(3a + 3a^q + 4)XY+
9e(a^q-a)(X+Y)+9e^2(a^q+a)=0.
\end{eqnarray*}
Note that if $3a + 3a^q + 4=0$ then $\deg(GCD(a^qX^3+X^2+b^q,bX^3+X+a))>0$, so we can assume now that $3a + 3a^q + 4\neq 0$.

Since $f_{a,b}(x)$ is a PP by hypothesis, then if $q\geq 47$ we have that $\mathcal{D}_{a,b}$ is reducible. Therefore the discriminant of $G_{a,b}^{\prime}(X,Y)$ with respect to $X$ is a square in $\overline{\mathbb{F}_{q}}[Y]$, which yields $3a^{q+1} + a + a^q=0$. Also, $-3$ cannot be a square in $\mathbb{F}_q$ otherwise the two linear components of $\mathcal{D}_{a,b}$ are $\mathbb{F}_q$-rational.

Summing up, we proved that if $f_{a,b}$ is a permutation and $\deg(\mathcal{D}_{a,b})<4$, then
\begin{equation}\label{SecondaTris}
a\ne -2/3, \qquad b=a^q+1/3,\qquad 3a^{q+1} + a + a^q=0.
\end{equation}
Finally, we prove that Condition \eqref{SecondaTris} implies Condition \eqref{SECONDA}.\\
Clearly, $b=-a^{q-1}/3$,  $b^{q+1}=1/9$, $a^q=-a/(3a+1)$.
Therefore, $$-3(1-4(b/a)^{q+1})=\frac{4-9a^{q+1}}{3a^{q+1}}=\frac{4+9a^2/(3a+1)}{-3a^2/(3a+1)}=-3\left(\frac{3a+2}{3a}\right)^2\in \mathbb{F}_q.$$
Since 

$$\left(-3\left(\frac{3a+2}{3a}\right)^2\right)^{(q-1)/2}=(-3)^{(q-1)/2}\left(\frac{3a+2}{3a}\right)^{q-1}=1,$$
Condition \eqref{SECONDA} holds.

\section{Appendix}\label{Sec:Appendix}

As already mentioned in the Introduction, another approach to establish whether a polynomial of type $f_{a,b}(X)$ permutes $\mathbb{F}_{q^2}$ is based on the classification of permutation rational functions of degree three on $\mathbb{P}_q^1=\mathbb{F}_q \cup \{\infty\}$ contained in \cite{FM2018}. Such a  technique, in principle, can be applied to determine if a polynomial $x^rh(X^{q-1})\in \mathbb{F}_{q^2}[X]$,  using (if any) classifications of rational functions of a fixed degree. For instance, in \cite{FM2018} the authors  provide a classification of rational functions of degree $3$ permuting $\mathbb{P}_q^1$ (up to composing on the left and on the right by M\"obius transformations). In particular, if $p>3$, there are only two classes: either $x^3$ or 
$$\varphi_{A}:=\frac{Ax^3+A^2x}{9x^2+A}, \quad    \textrm{ and } -A \textrm{ not a square in } \mathbb{F}_q;$$
see \cite[Table 1]{FM2018}. 

The whole class  of rational functions of degree $3$ permuting $\mu_{q+1}$ can be obtained by considering generic bijections from $\mu_{q+1}$ to $\mathbb{P}_q^1$ which have been classified by Zieve in \cite[Lemma 3.1]{Zieve2013}. In fact, a generic bijection $\ell(x): \mu_{q+1} \to \mathbb{P}_q^1$ can be written as
$$\ell_{\alpha,\beta}(x):=\frac{\alpha x-\beta \alpha^q}{x-\beta},$$
for some $\alpha \in \mathbb{F}_{q^2}\setminus \mathbb{F}_q$ and $\beta \in \mu_{q+1}$. More in details, a rational functions of degree $3$ permuting $\mu_{q+1}$ is of type

\begin{equation*}
\ell_{\gamma,\delta}^{-1} \circ \varphi_{A} \circ \ell_{\alpha,\beta}(x)=\ell_{\gamma,\delta}^{-1}\left(\varphi_{A}\left(\frac{\alpha x-\beta \alpha^q}{x-\beta}\right) \right)
\end{equation*}

or 

\begin{equation*}
\ell_{\gamma,\delta}^{-1} \left(\ell_{\alpha,\beta}(x)^3\right)=\ell_{\gamma,\delta}^{-1}\left(\left(\frac{\alpha x-\beta \alpha^q}{x-\beta}\right)^3 \right),
\end{equation*}

for some $\alpha,\beta,\gamma,\delta,A$. Let us consider only the first case. After easy computations, one gets 
\begin{equation}\label{Eq:Frac}
\frac{N_3 x^3+N_2x^2+N_1x+N_0}{D_3 x^3+D_2x^2+D_1x+D_0},
\end{equation}
where 
\begin{eqnarray*}
N_0&:=&-A^2\alpha^q\beta^3 \delta   - A \alpha^{3q}\beta^3 \delta  + A\beta^3 \delta  \gamma^q +
        9 \alpha^{2q} \beta^3 \delta \gamma^q;\\
N_1&:=&    A^2\alpha \beta^2 \delta  + 3A \alpha \beta^2 \delta  \alpha^{2q} - 
        18 \alpha^{q+1} \beta^2  \gamma^q\delta  + 2A^2\alpha^q \beta^2 \delta   - 
        3 A\beta^2   \gamma^q\delta - 9\alpha^{2q} \beta^2   \gamma^q\delta ;\\
N_2&:=&    -3A \alpha^{q+2} \beta \delta  + 9 \alpha^2 \beta  \gamma^q \delta- 
        2  A^2\alpha \beta \delta + 18 \alpha^{q+1} \beta \gamma^q\delta - A^2 
       \alpha^q \beta \delta + 3A  \beta \delta  \gamma^q;\\
N_3&:=&A \alpha^3 \delta - 9 \alpha^2  \gamma^q\delta +A^2 \alpha \delta  - A \gamma^q\delta ;\\
D_0&:=&    A\beta^3 \gamma  + 9 \alpha^{2q}\beta^3 \gamma  - A^2 \alpha^2q\beta^3  - 
        A \alpha^{3q}\beta^3 \\
D_1&:=&    -18 \alpha^{q+1} \beta^2 \gamma + A^2 \alpha \beta^2 + 3 A\alpha^{2q+1} \beta^2  
        - 3 A\beta^2 \gamma  - 9  \alpha^{2q}\beta^2 \gamma + 2 A^2 \alpha^q\beta^2 ;\\
D_2&:=&    9 \alpha^2 \beta \gamma - 3  A\alpha^{q+2} \beta + 18 \alpha^{q+1} \beta \gamma 
        - 2A^2 \alpha \beta  + 3A \beta \gamma  - A^2  \alpha^q\beta;\\
D_3&:=&    A \alpha^3 - 9 \alpha^2 \gamma +A^2 \alpha  -A \gamma.
\end{eqnarray*}
In general, it is not easy, for  fixed $a,b\in \mathbb{F}_{q^2}$, to determine  the values $A,\alpha,\beta,\gamma,\delta$ for which the rational function $\ell_{\gamma,\delta}^{-1} \circ \varphi_{A} \circ \ell_{\alpha,\beta}(x)$ equals $h_{a,b}(x)=\frac{a^qx^3+x^2+b^q}{bx^3+x+a}$. In fact, solving the system 
$$N_3^q=D_0, \quad N_2=1, \quad N_1=0, \quad N_0^q=D_3, \quad, D_2=0, \quad D_1=1$$
is not straightforward, since from a case by case analysis and eliminating $\beta,\gamma,\delta$ one obtains
\begin{eqnarray*}
9 \alpha^6 A^4 - 54 \alpha^{2q+6} A^3 + 81 \alpha^{4q+6} A^2 - 18 \alpha^{q+5} A^4 + 108 \alpha^{3q+5} A^3 - 162 \alpha^{5q+5} A^2 - 6 \alpha^4 A^5 \\
+ 45 \alpha^{2q+4} A^4 - 108 \alpha^{4q+4} A^3 + 81 \alpha^{6q+4} A^2  + 12 \alpha^{q+3} A^5 - 72 \alpha^{3q+3} A^4 + 108 \alpha^{5q+3} A^3+18 \alpha^{q+3}\\
 + \alpha^2 A^6 - 12 \alpha^{2q+2} A^5 + 45 \alpha^{4q+2} A^4  - 54 \alpha^{6q+2} A^3  +3 \alpha^2 A + 45 \alpha^{2q+2}- 2 \alpha^{q+1} A^6 + 12 \alpha^{3q+1} A^5 \\
- 18 \alpha^{5q+1} A^4  +12 \alpha^{q+1} A + 18 \alpha^{3q+1} + A^6 \alpha^{2q} - 6 A^5 \alpha^{4q} + 9 A^4 \alpha^{6q} + A^2 + 3 A\alpha^{2q}=0
\end{eqnarray*}
from which it is not easy to determine $A$ or $\alpha$. 

Also, once determined $a=D_0$ and $b=D_3$ in terms of $A,\alpha,\beta,\gamma,\delta$, one should still prove that Condition \eqref{PRIMA} or Condition \eqref{SECONDA} holds. This is the reason why, in order to determine the values $a,b\in \mathbb{F}_{q^2}$ for which $f_{a,b}(X)$ permutes $\mathbb{F}_{q^2}$, we used a different approach, based on algebraic curves.

\section{Acknowledgments*}
The research of D. Bartoli and M. Timpanella was partially supported  by the Italian National Group for Algebraic and Geometric Structures and their Applications (GNSAGA - INdAM).


\begin{thebibliography}{99}
\bibitem{Bartoli2018} D. Bartoli, On a conjecture about a class of permutation trinomials, Finite Fields Appl. {\bf 52}, 30--50 (2018).

\bibitem{BG2018} D. Bartoli, M. Giulietti, Permutation polynomials, fractional polynomials, and algebraic curves, Finite Fields Appl. {\bf 51}, 1 -- 16 (2018).

\bibitem{BQ2018} D. Bartoli, L. Quoos, Permutation polynomials of the type $x^rg(x^s)$ over $\mathbb{F}_{q^{2n}}$, Des. Codes Cryptogr. {\bf 86}, 1589--1599 (2018).

\bibitem{MAGMA}  W. Bosma, J. Cannon, C.  Playoust,  The Magma algebra system. {I}. The user language, J. Symbolic Comput. {\bf 24}, 235--265  (1997). 

\bibitem{DQWYY2015} C. Ding, L. Qu, Q. Wang, J. Yuan, P. Yuan, Permutation trinomials over finite fields with even characteristic, SIAM J. Discrete Math. {\bf 29}, 79--92 (2015).

\bibitem{FM2018} A. Ferraguti, G. Micheli, Full Classification of permutation rational functions and complete rational functions of degree three over finite fields, arxiv.org/abs/1805.03097 (2018)

\bibitem{GS2016} R. Gupta, R. K. Sharma, Some new classes of permutation trinomials over finite fields with even characteristic, Finite Fields Appl. {\bf 41}, 89 --96 (2016).

\bibitem{Hou2015} X. Hou, Determination of a type of permutation trinomials over finite fields, II, Finite Fields Appl. {\bf 35}, 16--35 (2015).

\bibitem{Hou2018} X. Hou, On a class of permutation trinomials in characteristic 2, Cryptogr. Commun. (2018). https://link.springer.com/article/10.1007/s12095-018-0342-1.

\bibitem{HTZ2019} X. Hou, Z. Tu, X. Zeng, Determination of a class of permutation trinomials in characteristic three, arXiv:1811.11949.

\bibitem{LP1997} J. B. Lee, Y. H. Park, Some permuting trinomials over finite fields, Acta Math. Sci. (English Ed.) {\bf 17}, 250--254 (1997).


\bibitem{LH2017} N. Li, T. Helleseth, Several classes of permutation trinomials from Niho exponents, Cryptogr. Commun. {\bf 9}, 693--705 (2017).

\bibitem{PL2001} Y. H. Park, J. B. Lee, Permutation polynomials and group permutation polynomials, Bull. Austral. Math. Soc. {\bf 63}, 67--74 (2001).

\bibitem{TZ2018} Z. Tu, X. Zeng, A class of permutation trinomials over finite fields of odd characteristic. Cryptography and Communications, available online May 1, 2018.

\bibitem{TZ2018b} Z. Tu, X. Zeng, Two classes of permutation trinomials with Niho exponents, Finite Fields Appl. {\bf 53}, 99--112 (2018).

\bibitem{TZLH2018}  Z. Tu, X. Zeng, C. Li, T. Helleseth, A class of new permutation trinomials. Finite Fields Appl. {\bf 50}, 178--195 (2018).

\bibitem{Sti} H. Stichtenoth,  Algebraic function fields and codes, Volume $254$ of Graduate Texts in Mathematics, 2nd edn. Springer, Berlin (2009).

\bibitem{Wang2007} Q. Wang, Cyclotomic mapping permutation polynomials over finite fields. In Sequences, Subsequences, and Consequences, S.W. Golomb, G. Gong, T. Helleseth, H.-Y. Song, (Eds.), pp. 119--128, Lecture Notes in Comput. Sci., vol. 4893, Springer, Berlin, 2007.

\bibitem{Zieve2009} M. E. Zieve, On some permutation polynomials over Fq of the form $x^rh(x(q-1)/d)$. Proc. Amer. Math. Soc. {\bf 137}, 2209--2216 (2009).

\bibitem{Zieve2013} M. E. Zieve, Permutation polynomials on $\mathbb{F}_q$ induced from bijective Redei functions on subgroups of the multiplicative group of $\mathbb{F}_q$, arxiv.org/abs/1310.0776 (2013).

\end{thebibliography}
\end{document}